\documentclass[a4paper,12pt]{article}

\usepackage{amsmath,amsfonts,amssymb,bm,amsthm}
\usepackage{mathrsfs}
\usepackage{cases}
\usepackage{cite}
\usepackage{indentfirst}
\usepackage[usenames]{color}
\usepackage{hyperref}
\usepackage{fancyhdr}
\usepackage{geometry}
\definecolor{webgreen}{rgb}{0,.5,0}
\definecolor{webbrown}{rgb}{.6,0,0}

\usepackage{color}

\newtheorem{theorem}{Theorem}
\newtheorem{lemma}{Lemma}

\newtheorem{definition}{Definition}

\newtheorem{proposition}{Proposition}
\newtheorem{corollary}{Corollary}

\allowdisplaybreaks

\begin{document}

\title{On the nonlinear Cauchy-Riemann equations of structural transformation and nonlinear Laplace equation}

\author{Gen Wang \thanks{{\tt
wanggen@zjnu.edu.cn}}  }

\date{\em \small {Department of Mathematics, Zhejiang Normal University, \\ ZheJiang, Jinhua 321004, P.R.China.} }

\maketitle

\begin{abstract}
This paper aims at studying a functional $K$-transformation  $w\left( z \right)\to \widetilde{w}\left( z \right)=w\left( z \right)K\left( z \right)$ that is made to reconsider the complex differentiability for a given complex function $w$ and subsequently we obtain structural holomorphic to judge a complex function to be complex structural differentiable. Since $K\left( z \right)$ can be chosen arbitrarily, thus it has greatly generalized the applied practicability.  And we particularly consider $K \left( z \right)= 1+\kappa \left( z \right)$,  then we found an unique Carleman-Bers-Vekua equations which is more simpler that all coefficients are dependent to the structural function $\kappa \left( z \right)$. The generalized exterior differential operator and the generalized Wirtinger derivatives are simultaneously obtained as well. As a discussion, second-order nonlinear Laplace equation is studied.

\vspace{.3 cm}
{{\bf{Keywords:}} Cauchy-Riemann equation, $K$-transformation, structural holomorphic condition, Carleman-Bers-Vekua equations, nonlinear Laplace equation }\\
{\bf{2010 Mathematics Subject Classification:}} 32V05, 32W05

\end{abstract}

\tableofcontents
\section{Introduction}

\subsection{Cauchy-Riemann equations}
In the field of complex analysis, the Cauchy-Riemann equations, consist of a system of two partial differential equations which, along with certain continuity and differentiability criteria, form a necessary and sufficient condition for a complex function to be complex differentiable, namely, holomorphic. This system of equations first appeared in the work of Jean le Rond d'Alembert in 1752. Later, Leonhard Euler in 1797 connected this system to the analytic functions. Cauchy in 1814 then used these equations to construct his theory of functions. Riemann's dissertation on the theory of functions appeared in 1851[see \cite{1,2,3,4,5}].

\begin{definition}[\cite{2}]\label{d2}
  Let $\Omega$ be an open set in $\mathbb {C}$ and $w$ a complex-valued function on $\Omega$. The
function $w$ is holomorphic at the point $z_{ 0} \in \Omega$ if the quotient
$\frac{w\left( {{z}_{0}}+h \right)-w\left( {{z}_{0}} \right)}{h}$
converges to a limit when $h\rightarrow 0$. Here $h\in \mathbb {C}$ and $h\neq 0$ with $z_{ 0} + h \in \Omega$, so that the quotient is well defined. The limit of the quotient, when it
exists, is denoted by $w'(z_{ 0} )$, and is called the derivative of $w$ at $z_{ 0 }$:
$$\lim _{{h\to 0}}{\frac {w(z_{0}+h)-w(z_{0})}{h}}=w'(z_{0})$$
\end{definition}
Emphasizedly, in the above limit, $h$ is a complex number
that may approach $0$ from any direction.
The function $w$ is said to be holomorphic on $\Omega$ if $w$ is holomorphic
at every point of $\Omega$. If $\Omega$ is a closed subset of $\mathbb {C}$, we say that $w$ is holomorphic on $\Omega$ if $f$ is holomorphic in some open set containing $\mathbb {C}$.
Finally, if $w$ is holomorphic in all of $\mathbb {C}$ we say that $w$ is entire.

Consider the complex plane $\mathbb {C} \equiv \mathbb {R} ^{2}=\{(x,y)\mid x\in \mathbb {R} ,\ y\in \mathbb {R} \}$. The Wirtinger derivatives are defined as the following linear partial differential operators of first order:
\[ {\frac {\partial }{\partial z}}={\frac {1}{2}}\left({\frac {\partial }{\partial x}}-\sqrt{-1}{\frac {\partial }{\partial y}}\right),\quad {\frac {\partial }{\partial {\bar {z}}}}={\frac {1}{2}}\left({\frac {\partial }{\partial x}}+\sqrt{-1}{\frac {\partial }{\partial y}}\right)\]
Clearly, the natural domain of definition of these partial differential operators is the space of $C^{1}$ functions on a domain $\Omega \subseteq \mathbb {R} ^{2}$, but, since these operators are linear and have constant coefficients, they can be readily extended to every space of generalized functions.

According to the definition \ref{d2},  the Cauchy-Riemann elliptic system of differential equations can be deducted as follows
\begin{equation}\label{eq21}
  \frac{\partial u}{\partial x}-\frac{\partial v}{\partial y}=0,~~\frac{\partial u}{\partial y}+\frac{\partial v}{\partial x}=0
\end{equation}
which is a classical method of construction of the theory of analytic functions $w=u+\sqrt{-1}v$ of a complex variable $z=x+\sqrt{-1}y$. For equation \eqref{eq21}, equivalently, the Wirtinger derivative of $w$ with respect to the complex conjugate of $z$ is zero,
\begin{equation}\label{eq19}
{\frac {\partial w}{\partial {\overline {z}}}}=0
\end{equation}
which is to say that, roughly, $w=u+\sqrt{-1}v$ is functionally independent from the complex conjugate of $z$.  The Cauchy-Riemann (CR) equations \eqref{eq21} plus equation \eqref{eq19} can be formally written as
\begin{equation}\label{q1}
  \left( \begin{matrix}
   \frac{\partial }{\partial x} & -\frac{\partial }{\partial y}  \\
   \frac{\partial }{\partial y} & \frac{\partial }{\partial x}  \\
\end{matrix} \right)\left( \begin{matrix}
   u  \\
   v  \\
\end{matrix} \right)=0\Leftrightarrow \overline{\partial }w=0
\end{equation}here $\overline{\partial }=d\overline{z}\frac{\partial }{\partial \overline{z}}$ is in $\mathbb {C}$.
The holomorphic functions coincide with those functions of two real variables with continuous first derivatives which solve the Cauchy-Riemann equations, a set of two partial differential equations.

\subsection{Carleman-Bers-Vekua (CBV) equation.}
The theory of generalized analytic functions founded by I. N. Vekua \cite{3} and L. Bers \cite{14} succeeded in being included in the pool of important techniques of the theory of partial differential equations. The reason is that the theory of generalized analytic functions is in a position to use the advantages of complex analysis for solving more general systems of partial
differential equations than this is possible in the framework of classical complex analysis.

Originally I. N. Vekua’s theory investigated only linear uniformly elliptic
systems for two desired real-valued functions in the plane. Today I. N. Vekua’s ideas are applied to partial differential equations in higher dimensions. Of course, a theory of the same high generality could not yet be developed so far\cite{8}.  Picard\cite{6}  proposed the idea of construction of a similar theory on the basis of a more general elliptic system of first order differential equations,
 \begin{align}\label{eq1}
  & {{\alpha }_{11}}\frac{\partial u}{\partial x}+{{\alpha }_{12}}\frac{\partial u}{\partial y}+{{\beta }_{11}}\frac{\partial v}{\partial x}+{{\beta }_{12}}\frac{\partial v}{\partial y}+{{a}_{1}}u+{{b}_{1}}v=0 \\
 & {{\alpha }_{21}}\frac{\partial u}{\partial x}+{{\alpha }_{22}}\frac{\partial u}{\partial y}+{{\beta }_{21}}\frac{\partial v}{\partial x}+{{\beta }_{22}}\frac{\partial v}{\partial y}+{{a}_{2}}u+{{b}_{2}}v=0 \notag
\end{align}
It is well known that under general assumptions about the coefficients, system  \eqref{eq1} is equivalent to the system
 \begin{equation}\label{eq17}
  \frac{\partial u}{\partial x}-\frac{\partial v}{\partial y}+au+bv=0,~~\frac{\partial u}{\partial y}+\frac{\partial v}{\partial x}+cu+dv=0
\end{equation}
 which was first investigated by Hilbert\cite{7} . Carleman \cite{9} obtained a fundamental property of the solutions of system \eqref{eq17}-their uniqueness. Earlier, Teodorescu \cite{10} studied a system of the following particular type:
 \[\frac{\partial u}{\partial x}-\frac{\partial v}{\partial y}+au+bv=0,~~\frac{\partial u}{\partial y}+\frac{\partial v}{\partial x}-bu+av=0\]
and obtained the general representation of the solutions by means of the analytic functions. This result turned out to be very important in constructing the general theory.

In the complex variables
 \begin{align}\label{eq22}
  & A=\frac{1}{4}\left( a+d+\sqrt{-1}c-\sqrt{-1}b \right),~~B=\frac{1}{4}\left( a-d+\sqrt{-1}c+\sqrt{-1}b \right)
\end{align}where unknown coefficients $a,b,c,d$ are chosen arbitrarily.
System \eqref{eq17} has the form
 \begin{equation}\label{eq20}
  \frac{\partial w}{\partial \overline{z}}+Aw+B\overline{w}=0
\end{equation}where  $w=u+\sqrt{-1}v$.
Equation \eqref{eq20} is called the Carleman-Bers-Vekua (CBV) equation.
 and its solutions are called generalized analytic functions; their theory is the meeting point between two sections of analysis-the theory of complex
variable analytic functions and the theory of elliptic type differential equations with two independent variables. The theory was developed as an independent part of analysis after appearance of the monograph of I. Vekua, where the long-term investigations of the author and some results of his disciples and followers (B. Bojarski, V. Vinogradov, Danilyuk, et al.) are presented. The foundations of the theory of generalized analytic functions were established in \cite{16}. Approximately in the same
period Bers, independently from Vekua, suggested a generalization of analytic functions (so-called pseudo-analytic functions), based on the modification of the concept of the derivative. Note that
many authors have proposed various generalizations, reducing system \eqref{eq1} to particular cases, until the complete theory of generalized analytic functions emerged[see \cite{20,21,22,24,25,26,27,28,29,30,31,32,33,34,35,36,37,38,39}].

Essentially, system \eqref{eq17} or equation \eqref{eq20} can be rewritten in a matrix form
\begin{equation}\label{eq8}
  \widehat{O}\left( \begin{matrix}
   u  \\
   v  \\
\end{matrix} \right)=0
\end{equation}
where $\widehat{O}=\left( \begin{matrix}
   \frac{\partial }{\partial x}+a & -\left( \frac{\partial }{\partial y}-b \right)  \\
   \frac{\partial }{\partial y}+c & \frac{\partial }{\partial x}+d  \\
\end{matrix} \right)$ is operator matrix.

A natural question arises while investigating the above-mentioned differential systems: how should the solutions be understood (the definition problem). It is clear that even for the simplest and most
fundamental case of system \eqref{eq21} it is not sufficient to assume the fulfillment of indicated differential equalities to obtain the class of functions with needed structure. For system \eqref{eq17} (Eq. \eqref{eq20}) the situation is even more complicated, since now additional coefficients are involved.  We will see below that this problem can be solved by $K$-transformation which makes it easier to be understood, especially, additional coefficients are endowed with specific meanings associated with one structural function $K$.

\subsection{the models of generalized Cauchy-Riemann system}
In the paper \cite{11},  the solvability of the Riemann-Hilbert problem for a generalized Cauchy-Riemann system with several singularities and reveal several new phenomenon has been studied by Heinrich Begehr  and Dao-Qing Dai,
they considered the generalized Cauchy-Riemann system
\[{{w}_{\overline{z}}}=\frac{Q\left( \overline{z} \right)}{P\left( \overline{z} \right)}w+aw+b\overline{w}\](see \cite{11} for more details). Let $G$ be a measurable set and $p$ be a real number, $1\le p\le \infty $. Denote by $L_{ p} (G)$ the set of all functions $w(z)$ satisfying the condition $\iint\limits_{G}{{{\left| w\left( z \right) \right|}^{p}}dxdy<+\infty }$, $L_{p}^{loc}\left( G \right)$ the set of all functions $w(z)$ satisfying the conditions ${{\left\| f \right\|}_{{{L}_{p}}\left( G \right)}}={{\left( \iint\limits_{G}{{{\left| w\left( z \right) \right|}^{p}}dxdy} \right)}^{\frac{1}{p}}}$, and $w\left( z \right)\in {{L}_{p}}\left( G' \right)$, respectively, where $G'$ is an arbitrary closed bounded subset of the set $G$.

For the first-order elliptic system of equations
\[{{u}_{\overline{z}}}=A\left( z \right)u+B\left( z \right)\overline{u}\]
with regular coefficients $A,B\in {{L}^{p}}\left( \Omega  \right),\left( p>2 \right)$  when $\Omega$ is a bounded domain in the complex plane $\mathbb{C}$,
 boundary value problems and their generalizations were investigated extensively over the past years,  This regularity allows to use a similarity principle for solutions of regular coefficient systems.

In shell theory, elliptic systems of equations with singular coefficients
occur
\[{{w}_{\overline{z}}}+\frac{A\left( z \right)}{\overline{z}}w+\frac{B\left( z \right)}{\overline{z}}\overline{w}=0\]
Now the coefficients do not belong to the regularity class $L^{ p}$ or $L^{ p;2}$ ; $p>2$: In \cite{12} a model equation ${{w}_{\overline{z}}}=\frac{b}{\overline{z}}\overline{w}$
was investigated, where $b$ is a complex constant. In \cite{13} the equation
${{w}_{\overline{z}}}+\frac{a\left( z \right)}{\left| z \right|}w+\frac{b\left( z \right)}{\left| z \right|}\overline{w}=0$
was studied. It was proved that there exist solutions admitting singularities of order $v>0$ at the point $z= 0$:  when $a=0$ was perturbed as
\[{{w}_{\overline{z}}}+\frac{b\left( 0 \right)}{\overline{z}}\overline{w}+\frac{b\left( z \right)-b\left( 0 \right)}{\overline{z}}\overline{w}=0\]
Under the assumption that $\frac{b\left( z \right)-b\left( 0 \right)}{\overline{z}}$ is sufficiently small, the existence of
continuous solutions was studied in \cite{12}.  In \cite{15,17}, \cite{11} found through the model equation
\[{{w}_{\overline{z}}}=\frac{\lambda }{\overline{z}}w+aw+b\overline{w};~~a,b\in {{L}^{p}}\left( \Omega  \right),~~p>2\]that the number of continuous solutions depends on size and sign of the constant $\lambda$:
This observation was implicitly supported by the results in \cite{18}, where the model equation
${{\psi }_{\overline{z}}}+\frac{a}{2\overline{z}}\psi +\frac{b}{2\overline{z}}\overline{\psi }=0$
was studied, where $a$ and $b$ are complex constants.   In \cite{19} the following boundary value problem was investigated:
\[\left\{ \begin{matrix}
   {{w}_{\overline{z}}}+A\left( z \right)w+B\left( z \right)\overline{w}=F,\left| z \right|<1  \\
   \operatorname{Re}\left[ {{z}^{-k}}w\left( z \right) \right]=g\left( z \right),\left| z \right|=1  \\
\end{matrix} \right.\]
[see more specific boundary value problems in \cite{11}].
the study about generalized Cauchy-Riemann system has brought so many significant results, but additional coefficients are still unknown and indeterminate.

\subsection{Nonlinear Cauchy-Riemann equations}
The Cauchy-Riemann equations are linear equations and they allow to solve only
linear Laplace equation. \cite{1} generalized C-R equations in such a
way to be able solve nonlinear Laplace equations such as system \eqref{eq17} or Carleman-Bers-Vekua (CBV) equation.

\begin{theorem}[NCR, \cite{1}]\label{t1}
If real functions $u(x,y)$ and $v(x,y)$ satisfy the NCR equations
\begin{equation}\label{eq4}
  {{u}_{y}}=-{{v}_{x}}+f\left( u,v \right),{{u}_{x}}={{v}_{y}}+g\left( u,v \right)
\end{equation}
where $g(u,v)$ and $f(u,v)$ are solutions of Cauchy-Riemann equations
\begin{equation}\label{eq3}
  {{f}_{u}}={{g}_{v}},{{f}_{v}}=-{{g}_{u}}
\end{equation}
then $u(x,y)$ and $v(x,y)$ satisfy nonlinear Laplace equations
\begin{align}\label{eq5}
  & \Delta u=\frac{1}{2}\frac{\partial }{\partial u}\left( {{f}^{2}}+{{g}^{2}} \right) \\
 & \Delta v=\frac{1}{2}\frac{\partial }{\partial v}\left( {{f}^{2}}+{{g}^{2}} \right) \notag
\end{align}

\end{theorem}
Let $u(x,y)$ and $v(x,y)$ are solutions of a system, which we call the nonlinear Cauchy-Riemann (NCR) equations, using \eqref{q1},  equation \eqref{eq4} can be shown in a matrix form
\[\left( \begin{matrix}
   \frac{\partial }{\partial x} & -\frac{\partial }{\partial y}  \\
   \frac{\partial }{\partial y} & \frac{\partial }{\partial x}  \\
\end{matrix} \right)\left( \begin{matrix}
   u  \\
   v  \\
\end{matrix} \right)=\left( \begin{matrix}
   g\left( u,v \right)  \\
   f\left( u,v \right)  \\
\end{matrix} \right)\]
where $f(u,v)$ and $g(u,v)$ are given functions.

\begin{lemma}[\cite{20}]\label{l1}
  The solution of the equation
  \begin{equation}\label{eq24}
    \frac{\partial w}{\partial \overline{z}}+Aw=0
  \end{equation}
on the whole plane, where $A\in L_{p}^{loc}\left( \mathbb{C} \right),p>2$,  has the form \[w\left( z \right)=\Phi \left( z \right){{e}^{-Q\left( z \right)}}\]where $Q(z)$ is one of $\frac{\partial }{\partial \overline{z}}$ -primitives of the function $A(z)$ and $\Phi(z)$ is an arbitrary entire function.
\end{lemma}
From the proof of lemma \ref{l1},
\[\frac{\partial w}{\partial \overline{z}}=-\Phi \left( z \right){{e}^{-Q\left( z \right)}}\frac{\partial Q}{\partial \overline{z}}=-w\left( z \right)\frac{\partial Q}{\partial \overline{z}}=-A\left( z \right)w\left( z \right)\]
it gives
$A\left( z \right)=\frac{\partial Q}{\partial \overline{z}}$,

\begin{theorem}[\cite{21}, the solution of 1-dim $\overline{\partial }$ problem]\label{t8}
  If $\eta \left( z \right)\in {{C}^{1}}\left( U \right)$ along with compact support, let
  \[h\left( z \right)=\frac{\sqrt{-1}}{2\pi }\iint\limits_{U}{\frac{\eta \left( \xi  \right)}{\xi -z}}d\overline{\xi }\wedge d\xi \]then $h\left( z \right)\in {{C}^{1}}\left( U \right)$, and it's a solution for $\frac{\partial h\left( z \right)}{\partial \overline{z}}=\eta \left( z \right)$.

\end{theorem}
The theorem \ref{t8}, traditionally, says that the inhomogeneous Cauchy-Riemann equations consist of the two equations for a pair of unknown functions $u(x,y)$ and $v(x,y)$ of two real variables
\[{\frac {\partial u}{\partial x}}-{\frac {\partial v}{\partial y}}=G (x,y),~~{\frac {\partial u}{\partial y}}+{\frac {\partial v}{\partial x}}=F (x,y)\]
for some given functions $G(x,y)$ and $F(x,y)$ defined in an open subset of $\mathbb{R}^{2}$. These equations are usually combined into a single equation
\begin{equation}\label{a2}
  {\frac {\partial w}{\partial {\bar {z}}}}=\varphi (z,{\bar {z}})
\end{equation}
where $w= u +\sqrt{-1}v$ and $\varphi = (G + \sqrt{-1}F)/2$.  If $\varphi$ is $\mathbb{C}^{k}$, then the inhomogeneous equation is explicitly solvable in any bounded domain $U$, provided $\varphi$ is continuous on the closure of $U$. Indeed, by the Cauchy integral formula,
$${\displaystyle w(\zeta ,{\bar {\zeta }})=\frac{\sqrt{-1}}{2\pi }\iint _{U}\varphi (z,{\bar {z}})\,{\frac {d{\bar {z}}\wedge dz}{z-\zeta }}}$$
for all $\zeta \in U$.

\section{Nonlinear Structural Cauchy-Riemann Equations}

\subsection{$K$-transformation}
For a given complex-valued function $w$ of a single complex variable, the derivative of $w$ at a point $z_{0}$ in its domain is defined by the limit \cite{1,2,3}  $ w'(z_{0})=\underset{z\to {{z}_{0}}}{\mathop{\lim }}{w(z)-w(z_{0}) \over z-z_{0}}$, if the limit exists, we say that $w$ is complex-differentiable at the point $z_{0}$.   Let's begin with the following definition.
\begin{definition}
Let $\Omega\subset \mathbb {C}$ be an open set and $w\left( z \right)=u+\sqrt{-1}v$ along with structural function $K\left( z \right)={{k}_{1}}+\sqrt{-1}{{k}_{2}}$ are complex valued function on $\Omega$, then a $K$-transformation such that
\begin{equation}\label{eq2}
  w\left( z \right)\to \widetilde{w}\left( z \right)=w\left( z \right)K\left( z \right)
\end{equation}
 where ${{k}_{1}},{{k}_{2}}$ are real functions with respect to the variables $x,y$, then
\begin{equation}\label{eq9}
  \widetilde{u}={{k}_{1}}u-v{{k}_{2}},~~\widetilde{v}=v{{k}_{1}}+u{{k}_{2}}
\end{equation}
hold in $\mathbb {C}$.
\end{definition}
If the real part takes the form ${{k}_{1}}=1+\alpha$, and ${{k}_{2}}=\beta $, then $K\left( z \right)=1+\kappa \left( z \right)$,
\[w\left( z \right)\to \widetilde{w}\left( z \right)=w\left( z \right)\left( 1+\kappa \left( z \right) \right)\]
where $\kappa \left( z \right)=\alpha +\sqrt{-1}\beta $ is structure function for all $z \in \Omega$ associated with the structure of specific manifold.

By the $K\left( z \right)=1+\kappa \left( z \right)$, the components of $ \widetilde{w}\left( z \right)= \widetilde{u}+\sqrt{-1}\widetilde{v}$  can be obtained
\begin{equation}\label{eq7}
  \widetilde{u}=u+u\alpha -v\beta ,~~\widetilde{v}=v+\alpha v+\beta u
\end{equation}
it also can be expressed in the matrix form
$$\left( \begin{matrix}
   \widetilde{u}  \\
   \widetilde{v}  \\
\end{matrix} \right)=\left( \begin{matrix}
   1+\alpha  & -\beta   \\
   \beta  & 1+\alpha   \\
\end{matrix} \right)\left( \begin{matrix}
   u  \\
   v  \\
\end{matrix} \right)$$Obviously, the functional structural matrix $\left( \begin{matrix}
   1+\alpha  & -\beta   \\
   \beta  & 1+\alpha   \\
\end{matrix} \right)$ is antisymmetric matrix. Therefore, the function
\begin{equation}\label{q11}
  \widetilde{w}\left( z \right)=\widetilde{u}+\sqrt{-1}\widetilde{v}=u+\alpha u-\beta v+\sqrt{-1}\left( v+\alpha v+\beta u \right)
\end{equation}
Note that the structural function is specified by specific manifolds.

{\bf{Remark 1}}:  As a convention, we always assume that real functions $u,v,\alpha ,\beta ,{{k}_{1}},{{k}_{2}}$ are continuous and differentiable for a better discussion. Analytic and holomorphic are treated as one explanation for complex functions. This paper will use Einstein summation convention $\sum\limits_{i}{d{{z}^{i}}\frac{\partial }{\partial {{z}^{i}}}}\equiv d{{z}^{i}}\frac{\partial }{\partial {{z}^{i}}}$.

\subsection{Structural holomorphic in  $\mathbb {C}$}
To derive the general expression relative to the  \eqref{eq2}, it leads to result similar to the general elliptic system of first order differential equations \eqref{eq1}
\begin{align}\label{eq18}
  & {{k}_{1}}\left( {{u}_{y}}+{{v}_{x}} \right)+{{k}_{2}}\left( {{u}_{x}}-{{v}_{y}} \right)+u\left( {{k}_{1y}}+{{k}_{2x}} \right)+v\left( {{k}_{1x}}-{{k}_{2y}} \right)=0 \\
 & {{k}_{1}}\left( {{u}_{x}}-{{v}_{y}} \right)-{{k}_{2}}\left( {{u}_{y}}+{{v}_{x}} \right)+u\left( {{k}_{1x}}-{{k}_{2y}} \right)-v\left( {{k}_{1y}}+{{k}_{2x}} \right)=0 \notag
\end{align}As mentioned, all coefficients are connected to the structural function $K\left( z \right)$, this is the most significant peculiarity of structural holomorphic. Equivalently, \eqref{eq18} can be expressed as
\[\left( \begin{matrix}
   {{k}_{1}} & {{k}_{2}}  \\
   -{{k}_{2}} & {{k}_{1}}  \\
\end{matrix} \right)\left( \begin{matrix}
   {{u}_{y}}+{{v}_{x}}  \\
   {{u}_{x}}-{{v}_{y}}  \\
\end{matrix} \right)=\left( \begin{matrix}
   -\left( {{k}_{1y}}+{{k}_{2x}} \right) & -\left( {{k}_{1x}}-{{k}_{2y}} \right)  \\
   -\left( {{k}_{1x}}-{{k}_{2y}} \right) & \left( {{k}_{1y}}+{{k}_{2x}} \right)  \\
\end{matrix} \right)\left( \begin{matrix}
   u  \\
   v  \\
\end{matrix} \right)\]or specifically
\[\left( \begin{matrix}
   {{k}_{1}} & {{k}_{2}}  \\
   -{{k}_{2}} & {{k}_{1}}  \\
\end{matrix} \right)\left( \begin{matrix}
   \frac{\partial }{\partial x} & \frac{\partial }{\partial y}  \\
   -\frac{\partial }{\partial y} & \frac{\partial }{\partial x}  \\
\end{matrix} \right)\left( \begin{matrix}
   v  \\
   u  \\
\end{matrix} \right)=\left( \begin{matrix}
   \frac{\partial }{\partial x} & \frac{\partial }{\partial y}  \\
   -\frac{\partial }{\partial y} & \frac{\partial }{\partial x}  \\
\end{matrix} \right)\left( \begin{matrix}
   -{{k}_{2}} & -{{k}_{1}}  \\
   -{{k}_{1}} & {{k}_{2}}  \\
\end{matrix} \right)\left( \begin{matrix}
   u  \\
   v  \\
\end{matrix} \right)\]
It is obvious that equation \eqref{eq18} is the same as system  \eqref{eq1} shown on the form, and equation \eqref{eq18} is equivalent to the equation
\begin{equation}\label{eq23}
  K\left( z \right)\frac{\partial w}{\partial \overline{z}}+w\frac{\partial K}{\partial \overline{z}}=0
\end{equation}
hence, if $K\left( z \right)=1+\kappa \left( z \right)=1+\alpha +\sqrt{-1}\beta $, then \eqref{eq23} is turned into \[\frac{\partial w}{\partial \overline{z}}+\kappa \frac{\partial w}{\partial \overline{z}}+w\frac{\partial \kappa }{\partial \overline{z}}=0\]
Therefore, the structural function $K$ has enormously opened our sights to explore various situations.

Putting equation \eqref{eq18} and system  \eqref{eq1} together is to derive the relation between coefficients
\begin{align}
  & {{\alpha }_{12}}={{\alpha }_{21}}={{\beta }_{11}}=-{{\beta }_{22}}={{k}_{1}} \notag\\
 & {{\beta }_{12}}={{\beta }_{21}}={{\alpha }_{22}}=-{{\alpha }_{11}}=-{{k}_{2}} \notag
\end{align}
and
\begin{align}
  & {{a}_{1}}=-{{b}_{2}}={{k}_{1y}}+{{k}_{2x}} \notag\\
 & {{b}_{1}}={{a}_{2}}={{k}_{1x}}-{{k}_{2y}} \notag
\end{align}
At this point, the structure function is $K={{k}_{1}}+\sqrt{-1}{{k}_{2}}$. Obviously, all the coefficients are related only to the component ${{k}_{1}},~{{k}_{2}}$ of the structure function $K$, and their first-order partial derivative, that is, to the structure function, which greatly simplifies the system \eqref{eq1}, and only needs to study the system \eqref{eq18}.

Note that $\kappa \frac{\partial w}{\partial \overline{z}}$ is very trivial to be considered,  in this paper, we don't take it into account, we need the nontrivial term $\frac{\partial \kappa }{\partial z},~\frac{\partial \kappa }{\partial \overline{z}}$ which are mainly subject to be studied by this paper.

\begin{corollary}
  The generalized structural Wirtinger derivatives are
  \begin{equation}\label{eq39}
  \frac{\rm{D}}{\partial z}=\frac{\partial }{\partial z}+\frac{\partial K}{\partial z},~~\frac{\rm{D}}{\partial \overline{z}}=\frac{\partial }{\partial \overline{z}}+\frac{\partial K}{\partial \overline{z}}
\end{equation}
  where $K\left( z \right)={{k}_{1}}+\sqrt{-1}{{k}_{2}}$ is complex function.

\end{corollary}

Generalized Carleman-Bers-Vekua equation can be formally written in the form
\begin{equation}\label{eq25}
  C\left( z \right)\frac{\partial w}{\partial \overline{z}}+wA\left( z \right)+B\left( z \right)\overline{w}=0
\end{equation}
Note that the parameters in the \eqref{eq25} are the same as \eqref{eq22} shown, the only difference is the function $C\left( z \right)$ added as a coefficient of first term, it is reduced to the Carleman-Bers-Vekua equation \eqref{eq20} if $C\left( z \right)=1$ holds.

\begin{definition}
  Let $\Omega\subset \mathbb {C}$ be an open set and complex valued function $w\left( z \right)=u+\sqrt{-1}v$ is said to be a structural holomorphic on $\Omega$ if and only if
\begin{equation}\label{eq36}
  \frac{\rm{D}}{\partial \overline{z}}w=0
\end{equation}
and its solutions are called generalized structural analytic functions.
\end{definition}
Let the collection of all structural holomorphic function be denoted by $Shol$, it means that if we say $w\in Shol$, then  $\frac{\rm{D}}{\partial \overline{z}}w=0$ holds.
Note that the $K\left( z \right)$ can be taken in different function form, by being taken in diverse form, it can suit distrinct needs.  Obviously, by lemma \ref{l1}, we can obtain the solution of \eqref{eq36} which is formally given by \[w\left( z \right)=\Phi \left( z \right){{e}^{-K\left( z \right)}}\]
Taking a expansion leads to
\[w\left( z \right)=\Phi \left( z \right)\left( 1-K\left( z \right)+\cdots  \right)=\Phi \left( z \right)+\Phi \left( z \right)h\left( z \right)\]
where $h\left( z \right)={{e}^{-K\left( z \right)}}-1$ is only induced by the structure function $K$.

\subsection{One case of structural holomorphic}

To start with complex differentiable in $\mathbb {C}$ based on the definition \ref{d2}. Firstly, we give a proposition based on the $K\left( z \right)=1+\kappa \left( z \right)$  to define the meaning of structural holomorphic.
\begin{proposition}\label{t4}
  Let $\Omega\subset \mathbb {C}$ be an open set and $ w\left( z \right)= u+\sqrt{-1}v$ a complex-valued continuous function on $\Omega$ is said to be structural holomorphic if
  \begin{equation}\label{eq32}
    \frac{{\rm{D}}w}{\partial z}=w'\left( z \right)+w\left( z \right)\kappa '\left( z \right)
  \end{equation}
 exists.
\begin{proof}
Given a complex valued function $w$ of a single complex variable, then based on the $K\left( z \right)=1+\kappa \left( z \right)$, accordingly, the derivative of $\widetilde{w}$ at a point $z_{0}$ in its domain is defined by the limit
\begin{equation}\label{eq6}
  \widetilde{w}'\left( {{z}_{0}} \right)=\underset{z\to {{z}_{0}}}{\mathop{\lim }}\,\frac{\widetilde{w}\left( z \right)-\widetilde{w}\left( {{z}_{0}} \right)}{z-{{z}_{0}}}
\end{equation}
Hence it gives
 $\widetilde{w}'\left( {{z}_{0}} \right)=\frac{{\rm{D}}w}{\partial z}\left( {{z}_{0}} \right)+\kappa \left( {{z}_{0}} \right)w'\left( {{z}_{0}} \right)$ for ${{z}_{0}}\ne z\in \Omega$, we say that $w$ is complex structural differentiable at the point $z_{0}$,  then
  $\widetilde{w}'\left( z \right)=\frac{{\rm{D}}w}{\partial z}+\kappa \left( z \right)w'\left( z \right)$, it directly indicates the existence of  $\frac{{\rm{D}}w}{\partial z}$.

\end{proof}
\end{proposition}

\begin{theorem}\label{t2}
Let $\Omega\subset \mathbb {C}$ be an open set and $w\left( z \right)=u+\sqrt{-1}v$ a complex-valued function on $\Omega$ is structural holomorphic if and only if
\begin{align}\label{eq15}
  &{{D}_{x}}u={{D}_{y}}v \\
  &{{D}_{y}}u=-{{D}_{x}}v \notag
\end{align}
where  ${{D}_{x}}=\frac{\partial }{\partial x}+{{\alpha }_{x}}-{{\beta }_{y}},~~{{D}_{y}}=\frac{\partial }{\partial y}+{{\alpha }_{y}}+{{\beta }_{x}}$.
\begin{proof}
As previously proved in the proposition \ref{t4}, complex structural differentiable \eqref{eq32} holds for all $z \in \Omega\subset \mathbb {C}$,
 then specifically
\begin{align}\label{eq13}
  & {{v}_{x}}+{{u}_{y}}+v\left( {{\alpha }_{x}}-{{\beta }_{y}} \right)+u\left( {{\beta }_{x}}+{{\alpha }_{y}} \right)=0 \\
 & {{u}_{x}}-{{v}_{y}}+u\left( {{\alpha }_{x}}-{{\beta }_{y}} \right)-v\left( {{\beta }_{x}}+{{\alpha }_{y}} \right)=0 \notag
\end{align}
It can be conveniently denoted as matrix form
$$\left( \begin{matrix}
   {{D}_{x}} & -{{D}_{y}}  \\
   {{D}_{y}} & {{D}_{x}}  \\
\end{matrix} \right)\left( \begin{matrix}
   u  \\
   v  \\
\end{matrix} \right)=0$$
where the generalized derivative operators are
\begin{equation}\label{q3}
  {{D}_{x}}=\frac{\partial }{\partial x}+{{\alpha }_{x}}-{{\beta }_{y}},~~{{D}_{y}}=\frac{\partial }{\partial y}+{{\alpha }_{y}}+{{\beta }_{x}}
\end{equation}
\end{proof}
\end{theorem}
In fact,
additional coefficients are involved of system \eqref{eq17} can be well interpreted by it, in comparison \eqref{eq8}, specifically,
\begin{align}\label{eq34}
  & a={{\alpha }_{x}}-{{\beta }_{y}},~~b=-\left( {{\alpha }_{y}}+{{\beta }_{x}} \right) \\
 & c={{\alpha }_{y}}+{{\beta }_{x}},~~d={{\alpha }_{x}}-{{\beta }_{y}} \notag
\end{align}
It turns out that identities $a=d,b=-c$ hold and are related to the structure function $\kappa$,  then it also leads to the further results of equation  \eqref{eq22}\[A=\frac{1}{2}\left( a+\sqrt{-1}c \right)=\frac{1}{2}\left( d-\sqrt{-1}b \right),~~~B=0\]then the Carleman-Vekua equation \eqref{eq20} is simplified as
$\frac{\partial w}{\partial \overline{z}}+Aw=0$.

Obviously, the equation \eqref{eq15} is the generalization of the C-R equation \eqref{q1}.
The nonlinear structural Cauchy-Riemann equations on a pair of real-valued functions of two real variables $u(x,y)$ and $v(x,y)$ are the two equations shown as equation \eqref{eq15}.  Thusly, based on the theorem \ref{t1}, we obtain
\begin{align}\label{a1}
  & f\left( u,v \right)=v\left( {{\beta }_{y}}-{{\alpha }_{x}} \right)-u\left( {{\beta }_{x}}+{{\alpha }_{y}} \right) \\
 & g\left( u,v \right)=v\left( {{\beta }_{x}}+{{\alpha }_{y}} \right)-u\left( {{\alpha }_{x}}-{{\beta }_{y}} \right) \notag
\end{align}and ${{f}_{v}}={{g}_{u}},~~{{f}_{u}}=-{{g}_{v}}$,
the generalized structural Wirtinger derivatives as an expansion of \eqref{eq8} can be correspondingly expressed as
\begin{equation}\label{eq10}
  \frac{{\rm{D}}}{\partial z}=\frac{1}{2}\left( {{D}_{x}}-\sqrt{-1}{{D}_{y}} \right),~~~ \frac{{\rm{D}}}{\partial \overline{z}} =\frac{1}{2}\left( {{D}_{x}}+\sqrt{-1}{{D}_{y}} \right)
\end{equation}
By simple calculation, \eqref{eq10} can properly be in a simpler form as corollary \eqref{eq39} shown,
\[\frac{{\rm{D}}}{\partial z}=\frac{\partial }{\partial {{z}}}+\frac{\partial \kappa }{\partial {{z}}},~~\frac{{\rm{D}}}{\partial \overline{z}}=\frac{\partial }{\partial {{\overline{z}}}}+\frac{\partial \kappa }{\partial {{\overline{z}}}}\]
then a simple theorem based on the theorem \ref{t2} can be naturally given in the following
\begin{corollary}\label{c1}
  Let $\Omega\subset \mathbb {C}$ be an open set and complex valued function $w\left( z \right)=u+\sqrt{-1}v$ is said to be a structural holomorphic on $\Omega$ if and only if
  \begin{equation}\label{eq16}
    \frac{{\rm{D}}}{\partial\overline{z}}w=0
  \end{equation}

\end{corollary}
Apparently, equation \eqref{eq16} as an expression is equivalent to equation \eqref{eq15}.  Traditionally, $\frac{\partial }{\partial \overline{z}}w=0$  is the fundamental foundation for the criterion to determine the analytic function such as $w$, known as C-R equation. In general, the C-R equation has its own restriction, in another words, it should be reasonably replaced by structural holomorphic equation. the operator $\overline{\partial }$ is accordingly replaced by the generalized $\overline{D }=\overline{\partial }+\overline{\partial} \kappa $. Importantly,  let's analyze the structural holomorphic condition of the corollary \ref{c1}, that is, the equation  $\frac{\partial w}{\partial \overline{z}}=-w\frac{\partial \kappa }{\partial \overline{z}}$.
By specific calculation, one obtains \[w\frac{\partial \kappa }{\partial \overline{z}}=\varphi \left( x,y \right)+\sqrt{-1}\zeta \left( x,y \right)\]where based on the equation \eqref{a1}, and
\begin{align}
  & \varphi \left( x,y \right)=\frac{1}{2}\left[ u\left( {{\alpha }_{x}}-{{\beta }_{y}} \right)-v\left( {{\alpha }_{y}}+{{\beta }_{x}} \right) \right]=-\frac{1}{2}G\left( u,v \right) \notag\\
 & \zeta \left( x,y \right)=\frac{1}{2}\left[ u\left( {{\alpha }_{y}}+{{\beta }_{x}} \right)+v\left( {{\alpha }_{x}}-{{\beta }_{y}} \right) \right]=-\frac{1}{2}F\left( u,v \right) \notag
\end{align}
it leads to the result as follows
\[\frac{\partial w}{\partial \overline{z}}=-\varphi \left( x,y \right)-\sqrt{-1}\zeta \left( x,y \right)\]or expressed as $\frac{\partial w}{\partial \overline{z}}=\frac{1}{2}\left( G+\sqrt{-1}F \right)$, this is equivalent to the equation \eqref{a2}, namely
$\varphi =\frac{1}{2}\left( G+\sqrt{-1}F \right)$, by the Cauchy integral formula, we get
\[w\left( \xi  \right)=\frac{\sqrt{-1}}{2\pi }\iint\limits_{U}{w(z)\frac{\partial \kappa }{\partial \overline{z}}}(z)\frac{dz\wedge d\overline{z}}{z-\xi }\]

On the foundation of lemma \ref{l1}, structural holomorphic in the corollary \ref{c1} has the similar expression as equation \eqref{eq24} shown,
$A=\frac{\partial \kappa }{\partial \overline{z}}$,
hence if $A=\frac{\partial \kappa }{\partial \overline{z}}\in L_{p}^{loc}\left( \mathbb{C} \right),p>2$ holds, then equation \eqref{eq16} has the solution
\begin{align}
w\left( z \right)  &=\Phi \left( z \right){{e}^{-\kappa \left( z \right)}}=\Phi \left( z \right){{e}^{-\alpha \left( z \right)}}{{e}^{-\sqrt{-1}\beta \left( z \right)}} \notag
\end{align}
In fact, one can choose the special expression to show the structure function
$\kappa \left( z \right)=\overline{z}A\left( z \right)$
for simple calculation $\frac{\partial \kappa }{\partial \bar{z}}=A\left( z \right)$.

\section{Functions of several complex variables}
There are Cauchy-Riemann equations, appropriately generalized, in the theory of several complex variables. They form a significant overdetermined system of PDEs. As often formulated, the d-bar operator  ${\bar {\partial }}$   annihilates holomorphic functions. This generalizes most directly the formulation   ${\partial w \over \partial {\bar {z}}}=0$, where ${\partial w \over \partial {\bar {z}}}={1 \over 2}\left({\partial w \over \partial x}+\sqrt{-1}{\partial w \over \partial y}\right)$.

Mathematically, the theory of functions of $n>1$ complex variables is to cope with complex valued functions $w=w(z^{1},z^2, \ldots, z^n)$ along with the structural function $\kappa=\kappa(z^{1},z^2, \ldots, z^n)$ on the space ${{\mathbb{C}}^{n}}$ of $n$-tuples of complex numbers for any $z={{\left( {{z}^{1}},{{z}^{2}},\cdots ,{{z}^{n}} \right)}^{T}}\in {{\mathbb{C}}^{n}}$, where the symbol $T$ represents the transpose of vectors. Definitely,  a structural holomorphic function can generalizes to several complex variables in a straightforward way. Let $\Omega\subset\mathbb {C} ^{n}$ denote an open subset, and let $w : \Omega \rightarrow \mathbb {C} $. Then function $w$ is analytic at a point $p$ in $\Omega$ if there exists an open neighborhood of $p$, where $w$ equals a convergent power series in $n$ complex variables, and meanwhile the nonlinear structural C-R equation still remains the fundamental importance on the complex manifold.

Consider the Euclidean space on the complex field $\mathbb {C} ^{n}=\mathbb {R} ^{2n}$, the Wirtinger derivatives are defined as the following matrix linear partial differential operators of first order:
\begin{equation}\label{a5}
\begin{split}
\frac{\partial }{\partial {{z}^{i}}}=\frac{1}{2}\left( \frac{\partial }{\partial {{x}^{i}}}-\sqrt{-1}\frac{\partial }{\partial {{y}^{i}}} \right)
\\ \frac{\partial }{\partial {{\overline{z}}^{i}}}=\frac{1}{2}\left( \frac{\partial }{\partial {{x}^{i}}}+\sqrt{-1}\frac{\partial }{\partial {{y}^{i}}} \right)
\end{split}
\end{equation}
As fundamental operator denoted, it's very useful to get the generalization further, hence to begin with the generalized Wirtinger derivatives.

For each index $i$ let $z^i=x^i+\sqrt{-1}y^i$, $w(z^i)=u^i+\sqrt{-1}v^i$, for each point $z=(z^1,\dots,z^n)\in U\subset\mathbb{C}^n$, and structural function takes the form $\kappa \left( {{z}^{i}} \right)=\alpha \left( {{x}^{i}},~~{{y}^{i}} \right)+\sqrt{-1}\beta \left( {{x}^{i}},{{y}^{i}} \right),~~i=1,\cdots ,n$,
and the generalized structural Cauchy-Riemann equation for one variable, then on the basis of generalized derivative operator \eqref{q3},  we obtain
\begin{corollary}\label{l2}
Generalized Wirtinger derivatives is  \[\frac{{\rm{D}}}{\partial z^{i}}=\frac{\partial }{\partial {{z}^{i}}}+\frac{\partial \kappa }{\partial {{z}^{i}}},~~\frac{{\rm{D}}}{\partial \overline{z^{i}}}=\frac{\partial }{\partial {{\overline{z}}^{i}}}+\frac{\partial \kappa }{\partial {{\overline{z}}^{i}}}\]
\end{corollary}
Note that this property implies that Wirtinger derivatives are derivations from the abstract algebra point of view, exactly unlike ordinary derivatives are, it's nonlinear because of the complex structural function $\kappa$.

\subsection{Structural holomorphic and $\overline{\rm{D}}$ operator }

On the space ${{\mathbb{C}}^{n}}$ of $n$-tuples of complex numbers for any $z={{\left( {{z}^{1}},{{z}^{2}},\cdots ,{{z}^{n}} \right)}^{T}}\in {{\mathbb{C}}^{n}}$. Similarly, a general $K$-transformation can be constructed as follows
\begin{equation}\label{eq26}
  w\left( z \right)\to \widetilde{w}\left( z \right)=w\left( z \right)K\left( z \right),~~ z\in U\subset {{\mathbb{C}}^{n}}
\end{equation}
where structural function $K\left( z \right)={{k}_{1}}+\sqrt{-1}{{k}_{2}}$, where ${{k}_{1}},{{k}_{2}}$ are real functions with respect to the variables $x=\left( {{x}^{1}},\cdots ,{{x}^{n}} \right),~~y=\left( {{y}^{1}},\cdots ,{{y}^{n}} \right)$, then
\begin{equation}
  \widetilde{u}={{k}_{1}}u-v{{k}_{2}},~~~\widetilde{v}=v{{k}_{1}}+u{{k}_{2}}
\end{equation}

\begin{corollary}\label{c2}
  The generalized Wirtinger derivatives based on \eqref{eq26} are
  \[\frac{{\rm{D}}}{\partial {z^{i}}}=\frac{\partial }{\partial z^{i}}+\frac{\partial K}{\partial z^{i}},~~\frac{{\rm{D}}}{\partial \overline{z^{i}}}=\frac{\partial }{\partial \overline{z^{i}}}+\frac{\partial K}{\partial \overline{z^{i}}}\]
  where $K\left( z \right)={{k}_{1}}+\sqrt{-1}{{k}_{2}}$ is complex function.

\end{corollary}

\begin{theorem}\label{t15}
  Let $\Omega\subset\mathbb{C}^n$ be an open set and complex valued function $w\left( z \right)=u+\sqrt{-1}v$ is said to be a structural holomorphic on $\Omega$ if and only if
  \begin{equation}\label{eq38}
    \frac{{\rm{D}}w}{\partial \overline{z^{i}}}=0
  \end{equation}
for all $z=(z^1,\dots,z^n)\in \Omega\subset\mathbb{C}^n$.

\end{theorem}
The theorem \ref{t15} implies that equation
$\frac{\partial w}{\partial {{{\bar{z}}}^{i}}}+w\frac{\partial \kappa }{\partial {{{\bar{z}}}^{i}}}=0,~~i=1,\cdots ,n$ holds for all structural holomorphic function $f$ on the $\Omega$, lemma \ref{l1} as an available solution can help it to be solved for each index $i$, its solution are $w\left( {{z}^{i}} \right)=\Phi \left( {{z}^{i}} \right){{e}^{-\kappa }}$.

Note that there is a special case, if let $w=K$ be given and taking it to the
theorem \ref{t15}, it yields $\frac{{\rm{D}}K}{\partial \overline{z^{i}}}=0$, and then it consequently deduces
\begin{equation}\label{eq28}
  \frac{\partial K}{\partial \overline{{{z}^{i}}}}=0
\end{equation}
as such special case above, the Wirtinger derivative of $K$ with respect to the complex conjugate of $z$ is zero, this is a very interesting feature of the structural function $K\left( z \right)$.

Actually, the theorem \ref{t15} is a unification which can work for complex field $\mathbb {C}$ or $\mathbb {C}^{n}$ or complex manifold $M$, if $n=1$, then $\mathbb {C}^{1}=\mathbb {C}$.  Thinking about structural holomorphic equation \eqref{eq38} \[\frac{\partial w}{\partial \overline{{{z}^{i}}}}=-w\frac{\partial K}{\partial \overline{{{z}^{i}}}}\]
if $\frac{\partial K}{\partial \overline{{{z}^{i}}}}=0$, then $\frac{\partial w}{\partial \overline{{{z}^{i}}}}=0$, it is rightly just about classical Cauchy-Riemann equation.

The generalized exterior differential form is shown as below theorem
\begin{theorem}\label{t7}
The generalized exterior differential operator is
\[\mathcal{D}=d+dK=\rm{D}+\overline{\rm{D}}\]where ${\rm{D}}=\partial +\partial K,~~\overline{\rm{D}}=\overline{\partial }+\overline{\partial }K$, and $d=\partial +\overline{\partial }$.
\begin{proof}The generalized exterior differential of complex valued functions $w$ is given by
\begin{align}
 \mathcal{D}w & =\frac{{\rm{D}}w}{\partial {z^{i}}}d{{z}^{i}}+\frac{{\rm{D}}w}{\partial \overline{z^{i}}}d\overline{{{z}^{i}}} \\
 & =\frac{\partial w}{\partial {{z}^{i}}}d{{z}^{i}}+w\frac{\partial K}{\partial {{z}^{i}}}d{{z}^{i}}+\frac{\partial w}{\partial \overline{{{z}^{i}}}}d\overline{{{z}^{i}}}+w\frac{\partial K}{\partial \overline{{{z}^{i}}}}d\overline{{{z}^{i}}} \notag\\
 & =\left( \frac{\partial w}{\partial {{z}^{i}}}d{{z}^{i}}+\frac{\partial w}{\partial \overline{{{z}^{i}}}}d\overline{{{z}^{i}}} \right)+w\left( \frac{\partial K}{\partial {{z}^{i}}}d{{z}^{i}}+\frac{\partial K}{\partial \overline{{{z}^{i}}}}d\overline{{{z}^{i}}} \right)\notag \\
 &=\left( \partial w+\overline{\partial }w \right)+w\left( \partial K+\overline{\partial }K \right)\notag \\
 & =dw+wdK \notag
\end{align}
where $dw=\frac{\partial w}{\partial {{z}^{i}}}d{{z}^{i}}+\frac{\partial w}{\partial \overline{{{z}^{i}}}}d\overline{{{z}^{i}}}=\partial w+\overline{\partial }w$. Then
 \[\mathcal{D}=d+dK=\left( \partial +\partial K \right)+\left( \overline{\partial }+\overline{\partial }K \right)=\rm{D}+\overline{\rm{D}}\]
 where ${\rm{D}}=\partial +\partial K,~~\overline{\rm{D}}=\overline{\partial }+\overline{\partial }K$.

\end{proof}
\end{theorem}

Basically, the theorem \ref{t15} can be viewed in a equivalent expression by using $\overline{\rm{D}}$ operator,

\begin{corollary}\label{c4}
The complex function  $w\left( z \right)=u+\sqrt{-1}v$ is structural holomorphic on $\Omega\subset \mathbb {C}^{n}$ is
\[\overline{\rm{D}}w=0\]
\end{corollary}

To see how the theorem \ref{t15} or corollary \ref{c4} unifies all situations, let's take the theorem \ref{t15} into consideration. Given a complex function $w$ defined on complex field $\Omega$, then a category can be created based on the structural function $K(z)$ as follows
\begin{description}
  \item[i] when $K(z)=1$ is constant number, then corollary \ref{c2} is rewritten as $\frac{{\rm{D}}}{\partial {z^{i}}}=\frac{\partial }{\partial {{z}^{i}}},\frac{{\rm{D}}}{\partial \overline{z^{i}}}=\frac{\partial }{\partial \overline{{{z}^{i}}}}$, this is classic representation shown as \eqref{a5}, meanwhile, structural holomorphic condition is reduced to classic Cauchy-Rieman equation $\frac{\partial w}{\partial \overline{{{z}^{i}}}}=0$. This is for the  $\mathbb {C}^{n}$, when it comes to $n=1$, it's also applied to the  $\mathbb {C}$, at this time, the operator matrix $\widehat{O}$ in \eqref{q1} is $\hat{O}=\left( \begin{matrix}
   \frac{\partial }{\partial x} & -\frac{\partial }{\partial y}  \\
   \frac{\partial }{\partial y} & \frac{\partial }{\partial x}  \\
\end{matrix} \right)$. Perspectively, in a viewpoint of structural holomorphic, the classical Cauchy-Rieman equation is equivalent to constant structural holomorphic condition.

  \item[ii] when $K\left( z \right)=1+\kappa \left( z \right)$ holds, then corollary \ref{c2} becomes corollary \ref{l2}.
Accordingly,  theorem \ref{t15} turns to theorem \ref{t15}. The operator matrix $\widehat{O}$ in \eqref{eq8} becomes antisymmetric operator matrix
  $\hat{O}=\left( \begin{matrix}
   \frac{\partial }{\partial x}+a & -\left( \frac{\partial }{\partial y}+c \right)  \\
   \frac{\partial }{\partial y}+c & \frac{\partial }{\partial x}+a  \\
\end{matrix} \right)$, where coefficients $a,c$ are as \eqref{eq34} shows.

\end{description}
Note that by introducing the structural function $K(z)$, on the one hand, it can unify all possible kind of equations to describe the holomorphic or generalized analytic, and it suits the complex manifold even higher dimensional manifold, on the other hand, we can study different kinds of situations as we need owing to structural function $K(z)$ which can be chosen arbitrarily. In other words, we can deduce all equations like Cauchy-Riemann equation or nonlinear Cauchy-Riemann equation or Carleman-Bers-Vekua equation by choosing appropriate structural function $K(z)$ defined on some complex domain.

The most important thing is that we notice that by using structural function $K(z)$, we can nicely study some vital properties about singularity or any special points in complex field.\\
\subsection{Examples}
{\bf{Example 1.}} Consider the structural function $K\left( z \right)={{e}^{z\overline{z}}}={{e}^{{{\left| z \right|}^{2}}}}$ in $\mathbb {C}$, then the generalized structural Wirtinger derivatives $\frac{{\rm{D}}}{\partial \overline{z}}$ can be shown as
\[\frac{{\rm{D}}}{\partial \overline{z}}=\frac{\partial }{\partial \overline{z}}+\frac{\partial K}{\partial \overline{z}}=\frac{\partial }{\partial \overline{z}}+z{{e}^{{{\left| z \right|}^{2}}}}\]
then structural holomorphic condition is \[\frac{\partial w}{\partial \overline{z}}+zKw=0\]Consider point $z=0$, then $\frac{\partial w}{\partial \overline{z}}=0$ holds at the origin.\\

{\bf{Example 2.}} If we consider $K\left( z \right)=\overline{z}$, then structural holomorphic condition on the whole plane can be represented in the following
form
\[\frac{{\rm{D}}w}{\partial \overline{z}}=\frac{\partial w}{\partial \overline{z}}+w=0\]Thusly, then $\frac{\partial w}{\partial \overline{z}}=-w$.\\

{\bf{Example 3.}} If we consider $K\left( z \right)={{e}^{{{\left| z \right|}^{2}}}}+\overline{z}$, then structural holomorphic condition on the whole plane is
\[\frac{{\rm{D}}w}{\partial \overline{z}}=\frac{\partial w}{\partial \overline{z}}+wz{{e}^{{{\left| z \right|}^{2}}}}+w=0\]
in this way, we have
$\frac{\partial w}{\partial \overline{z}}=-wz{{e}^{{{\left| z \right|}^{2}}}}-w$, equation $\left( \frac{\partial w}{\partial \overline{z}}+w \right)\left| _{z=0} \right.=0$ holds at origin $z=0$.\\

\section{ Nonlinear Laplace equation}
\begin{proposition}
  Nonlinear structural Laplace operator on $\Omega\subset \mathbb {C}^{n}$ is
  \begin{equation}\label{eq27}
    \frac{{{\rm{D}}^{2}}}{\partial{{z}^{i}}\partial\overline{{{z}^{j}}}}=\frac{{{\partial }^{2}}}{\partial {{z}^{i}}\partial \overline{{{z}^{j}}}}+\left( \frac{\partial K}{\partial {{z}^{i}}}\frac{\partial }{\partial \overline{{{z}^{j}}}}+\frac{\partial K}{\partial \overline{{{z}^{j}}}}\frac{\partial }{\partial {{z}^{i}}} \right)+\left( \frac{{{\partial }^{2}}K}{\partial {{z}^{i}}\partial \overline{{{z}^{j}}}}+\frac{\partial K}{\partial {{z}^{i}}}\frac{\partial K}{\partial \overline{{{z}^{j}}}} \right)
  \end{equation}
\begin{proof}
  For a given function  $w\left( z \right)=u+\sqrt{-1}v$ that is defined on the $\Omega\subset \mathbb {C}^{n}$, then based on the \ref{c2}, one obtains
\[\frac{{{\rm{D}}^{2}}}{\partial{{z}^{i}}\partial\overline{{{z}^{j}}}}=\frac{{{\partial }^{2}}}{\partial {{z}^{i}}\partial \overline{{{z}^{j}}}}+\left( \frac{\partial K}{\partial {{z}^{i}}}\frac{\partial }{\partial \overline{{{z}^{j}}}}+\frac{\partial K}{\partial \overline{{{z}^{j}}}}\frac{\partial }{\partial {{z}^{i}}} \right)+\left( \frac{{{\partial }^{2}}K}{\partial {{z}^{i}}\partial \overline{{{z}^{j}}}}+\frac{\partial K}{\partial {{z}^{i}}}\frac{\partial K}{\partial \overline{{{z}^{j}}}} \right)\]
Therefore, one gets the desire result.
\end{proof}
\end{proposition}
One notices that nonlinear Laplace operator is dependent to the specific form of the generalized Wirtinger derivatives in corollary \ref{c2}.

Take notice of something that nonlinear Laplace operator is partial differential operator of second order, it well fits almost all situation in complex field. It appears that it's an unification for all possible situation.
Obviously, if let $K=1$ be given, then nonlinear Laplace operator on $\Omega\subset \mathbb {C}^{n}$ will be degenerately rewritten as
\[\frac{{{\rm{D}}^{2}}}{\partial{{z}^{i}}\partial\overline{{{z}^{j}}}}  =\frac{{{\partial }^{2}}}{\partial {{z}^{i}}\partial \overline{{{z}^{j}}}}\]it's corresponding to the classical Laplace operator.

Denote ${{K}_{i\overline{j}}}=\frac{{{\partial }^{2}}K}{\partial {{z}^{i}}\partial \overline{{{z}^{j}}}}$, ${{\psi }_{i\overline{j}}}={{K}_{i\overline{j}}}+\frac{\partial K}{\partial {{z}^{i}}}\frac{\partial K}{\partial \overline{{{z}^{j}}}}$, then \eqref{eq27} can be rewritten in the form
\begin{equation}\label{eq30}
  \frac{{{\rm{D}}^{2}}}{\partial{{z}^{i}}\partial\overline{{{z}^{j}}}}=\frac{{{\partial }^{2}}}{\partial {{z}^{i}}\partial \overline{{{z}^{j}}}}+\frac{\partial K}{\partial \overline{{{z}^{j}}}}\frac{\partial }{\partial {{z}^{i}}}+\frac{\partial K}{\partial {{z}^{i}}}\frac{\partial }{\partial \overline{{{z}^{j}}}}+{{\psi }_{i\overline{j}}}
\end{equation}
As previously demonstrated, one easily Laplace equations $\frac{{{\partial}^{2}}w}{\partial{{z}^{i}}\partial\overline{{{z}^{j}}}}=0$, that is,
\begin{theorem}\label{t6}
  For structural holomorphic function $w$ is defined on the $\Omega\subset \mathbb {C}^{n}$ such that
\begin{equation}\label{eq29}
  \frac{{{\rm{D}}^{2}}}{\partial{{z}^{i}}\partial\overline{{{z}^{j}}}}=0
\end{equation}is called nonlinear Laplace equations.
\end{theorem}
The second-order partial differential equations \eqref{eq29} is a differential equation that contains unknown multivariable functions $K$ and their partial derivatives $\frac{\partial K}{\partial \overline{{{z}^{j}}}},\frac{\partial K}{\partial {{z}^{i}}}$. It's obvious that it has different types in different regions with different structural function $K$. The structural function $K$ provides a guide to appropriate initial and boundary conditions, and to the smoothness of the solutions. There are examples of linear partial differential equations whose coefficients have derivatives of all orders but which have no solutions at all. Even if the solution of a partial differential equation exists and is unique, it may nevertheless have undesirable properties.

Equations \eqref{eq29} as a PDEs actually can be used to describe a wide variety of phenomena such as sound, heat, electrostatics, electrodynamics, fluid dynamics, elasticity, or quantum mechanics by choosing proper structural function $K(z)$ on the region we consider. These seemingly distinct physical phenomena can be formalised similarly in terms of PDEs. .

Transparently, \eqref{eq29} is partial differential equation of second order with form invariance of equation.  By using the operators $\partial ,\overline{\partial }$ and exterior differential operator $\wedge$, equation \eqref{eq29} can be reshown as
\begin{align}\label{eq31}
  \frac{{{\rm{D}}^{2}}}{\partial{{z}^{i}}\partial\overline{{{z}^{j}}}}d{{z}^{i}}\wedge d\overline{{{z}^{j}}}
 &=\partial \overline{\partial }w+\partial w\overline{\partial }K+\partial K\overline{\partial } w+w\psi\\
 & =0 \notag
\end{align}where \[\psi ={{\psi }_{i\overline{j}}}d{{z}^{i}}\wedge d\overline{{{z}^{j}}}=\partial \overline{\partial }K+\partial K\overline{\partial }K\]
As shown in the equation \eqref{eq30} or \eqref{eq31}  , the last term ${{\psi }_{i\overline{j}}}$ or $\psi$ is a sort of metric function, not a operator, because of this unique feature, it's worth deeply analyzing it. By the way, we notice that theorem \ref{t6} is a form of second order, it reveals that   \eqref{eq29} is a self-contained equation, it means equation \eqref{eq30} is an invariant second order operator under the $K$-transformation between complex function on complex field $\mathbb {C}$ or $\mathbb {C}^{n}$ or complex manifold $M$.

For the functions of one complex variable, and based on the theorem \ref{t6}, it directly derives the nonlinear Laplace equations with one complex variable. To begin with a nonlinear Laplace operator
\[{{\Delta }_{K}}=\frac{{1}}{4}\Delta +\frac{\partial K}{\partial \overline{z}}\frac{\partial }{\partial z}+\frac{\partial K}{\partial z}\frac{\partial}{\partial \overline{z}}+\eta \]where $\Delta =4\frac{{{\partial }^{2}}}{\partial z\partial \overline{z}}$ is Laplace operator, $\eta =\frac{1}{4}\Delta K+\frac{\partial K}{\partial z}\frac{\partial K}{\partial \overline{z}}$.

\section{Conclusions}
Taking it by and large, this paper has successfully built a clear theory mode to unify all kinds of theoretic form by doing the transformation, actually, it's a unified transformation. Mathematically, $\kappa$-transformation is a special case of $K$-transformation. Subsequently, nonlinear Laplace equation is developed based on the generalized Wirtinger derivatives and structural holomorphic condition. Hence we will separately explain in the following

\begin{description}
  \item[$K=1$] Cauchy-Riemann equation, analytic condition in $\mathbb {C}^{n}$ \[\frac{\partial w}{\partial \overline{{{z}^{i}}}}=0,\left( i=1,\cdots ,n \right)\]that is, structural holomorphic condition.
  \item[$K=1+\kappa$]  nonlinear structural Cauchy-Riemann equations, generalized analytic condition in $\mathbb {C}^{n}$\[\frac{\partial w}{\partial \overline{{{z}^{i}}}}+w\frac{\partial \kappa }{\partial \overline{{{z}^{i}}}}=0,\left( i=1,\cdots ,n \right)\]
   that is, structural holomorphic condition.

\end{description}
In $\mathbb {C}$, $n=1$ is taken and the equation form remains.

\end{document}